\documentclass[12pt]{article}
\usepackage{amsmath}
\usepackage{amsfonts}
\usepackage{amsgen}
\usepackage{amstext}
\usepackage{amssymb}

\begin{document}
\setlength{\baselineskip}{5mm}
\setlength{\topmargin}{1.2cm}
\setlength{\textheight}{19.0cm}
\setlength{\oddsidemargin}{2.4cm}
\setlength{\evensidemargin}{2.4cm}
\setlength{\textwidth}{11.5cm}
\pagestyle{headings}
\newtheorem{theorem}{Theorem}
\newtheorem{corollary}[theorem]{Corollary}
\newtheorem{example}[theorem]{Example}
\newtheorem{proposition}[theorem]{Proposition}
\newtheorem{remark}[theorem]{Remark} 
\newtheorem{definition}[theorem]{Definition} 
\newtheorem{lemma}[theorem]{Lemma}
\newtheorem{conjecture}[theorem]{Conjecture}
\newtheorem{algorithm}[theorem]{Algorithm}
\renewcommand{\thetheorem}{\arabic{theorem}}
\def\ppp{{\mathbb{P}}}
\def\aaa{{\mathbb{A}}}
\def\fff{{\mathbb{F}}}
\def\qqq{\mathbb{Q}}
\def\rrr{\mathbb{R}}
\def\ccc{\mathbb{C}}
\def\hhh{\mathbb{H}}
\def\zzz{\mathbb{Z}}
\def\ddd{\,}
\def\nnn{\mathbb{N}}
\def\ggg{\mathbb{G}}
\def\pf{{\bf proof}:\ }
\def\qed{$\Box$}
\def \dis {\displaystyle}
\def \dotsc {\dots}

\author{David Joyner and Salahoddin Shokranian}
\title{Remarks on codes from modular curves: MAGMA applications\thanks{This
paper is a modification of the paper \cite{JS}. Here we use 
\cite{MAGMA}, version 2.10, instead of MAPLE. Some examples have been changed
and minor corrections made.
MSC 2000: 14H37, 94B27,20C20,11T71,14G50,05E20,14Q05 }}
\date{3-30-2004}
\maketitle

\tableofcontents

\section{Introduction}

        Suppose that $V$ is a smooth projective
variety over a finite field $k$.
An important problem in arithmetical algebraic geometry
is the calculation of the number of $k$-rational
points of $V$, $|V(k)|$. The work of Goppa \cite{G}
and others have shown its importance
in geometric coding theory as well. 
We refer to this problem as the {\bf counting problem}. In
\index{counting problem}
most cases it is very hard to find an explicit formula
for the number of points of
a variety over a finite field.

When the variety is a
``Shimura variety'' defined by certain
group theoretical conditions (see \S \ref{section:shimura} below),
methods from non-abelian harmonic
analysis on groups can be used to find an explicit solution for the
counting problem. The Arthur-Selberg trace formula \cite{S}, provides one
such method.
Using the Arthur-Selberg trace formula, an
explicit formula for the counting problem has been found
for Shimura varieties, thanks primarily to the work of Langlands and Kottwitz
(\cite{Lan1},  \cite{K1})
\footnote{For some introductions to this highly
technical work of Langlands and Kottwitz, the reader is referred
to Labesse \cite{Lab}, Clozel \cite{Cl}, and
Casselman \cite{Cas2}.}.
Though it may be surprising and indeed
very interesting that the trace formula allows one
(with sufficient skill and expertise) to relate, when $V$ is
a Shimura variety, the geometric numbers $|V({\fff}_q)|$ to
orbital integrals from harmonic analysis (\cite{Lab}, for example),
or to a linear combination of coefficients
of automorphic forms (\cite{Gel}, for example), or even to
representation-theoretic data (\cite{Cas2}, for example),
these formulas do not yet seem to be helping the coding theorist
in any practical way that we know of.

However, another type of application of the trace formula
is very useful. Moreno \cite{M}
first applied the trace formula in the context of Goppa codes
to obtaining a new proof of a famous result of
M. Tsfasman , S. Vladut, T. Zink, and Y. Ihara. (Actually, Moreno
used a formula for the trace of the Hecke operators
acting on the space of modular forms of weight 2, but this can
be proven as a consequence of the Arthur-Selberg trace formula,
\cite{DL}, \S II.6.) This will be discussed below.
We are going to restrict our attention in this paper to the interplay
between Goppa codes of modular curves and
the counting problem, and give some examples using MAGMA.
In coding theory, curves with many rational points over
finite fields are being used for construction of codes with some good
specific characteristics. We discuss AG (or Goppa) codes arising from curves,
first from an abstract general perspective then turning
to concrete examples associated to modular curves.
We will try to explain these extremely technical ideas using a special case
at a level to a typical graduate student with some background
in modular forms, number theory, group theory, and algebraic geometry.
For an approach similar in spirit, though from a more classical
perspective, see the book of C. Moreno \cite{M}.

\section{Shimura curves}
\label{section:shimura}

        In this section we study arithmetic subgroups, arithmetical
quotients, and their rational compactifications. Ihara first introduced
Shimura curves, a rational compactification of $\Gamma\backslash \hhh$
where
$\Gamma$ is a particular discrete subgroup, from a classical perspective.
We shall recall them from both the classical and group-theoretical point 
of
view. The latter perspective generalizes to higher dimensional Shimura
varieties \cite{Del}.

\subsection{Arithmetic subgroups}

We assume that $G = SL(2)$ is the  group of $2 \times 2$ matrices with
entries from an algebraically closed field $\Omega$. In particular the
group of $R$-points of $SL(2)$ for a subring $R \subseteq \Omega$, with
unit element $1$ is defined by
\[
SL(2, R) =\{ g \in M(2, R)\;|\; \det(g) =1 \},
\]
where $M(2, R)$ is the space of $2 \times 2$ matrices with entries from
$R$. We now define congruence subgroups in $SL(2, \zzz)$. Let $SL(2,\zzz)$
be the subgroup of $SL(2, \rrr)$ with integral matrices. Consider a
natural number $N$, and let
\[
\Gamma(N) = \left\{ \left[\begin{array}{cc}
a & b\\
c & d
\end{array} \right] \in SL(2, \zzz)\;|\; \begin{array}{rr}
a, d \equiv 1 ({\rm mod}\, N)\\
b , c \equiv 0 ({\rm mod}\, N)
\end{array}
\right\},
\]
We note that the  subgroup $\Gamma(N)$ is  a discrete subgroup of
$SL(2, \rrr)$, which is called the {\em principal congruence subgroup of
level $N$}. Any subgroup of $SL(2, \zzz)$ that contains the principal
congruence subgroup is called a {\em congruence subgroup}.

        In general an {\em arithmetic subgroup} of $SL(2, \rrr)$ is any
discrete subgroup $\Gamma$ that is commensurable with $SL(2, \zzz)$, where
{\em commensurability} means that the intersection $\Gamma \cap SL(2,\zzz)$
is of finite index in both $\Gamma$ and $SL(2, \zzz)$. The group
$\Gamma(N)$ has the property of being commensurable with
$SL(2, \zzz)$.

\subsection {Riemann surfaces as algebraic curves}

        Let us recall that the  space  $\hhh =\{ z \in \ccc\;|\; Im(z) > 0\}$
is called {\em the Poincar\'e upper half plane}. This space plays
fundamental r\^ole in the definition of the modular curves. Note that the
group $SL(2, \rrr)$ acts on $\hhh$ by
\[
g\cdot z = (az + b) (cz + d)^{-1} = \frac{az + b}{cz +d},
\]
where $z \in \hhh$, $g = \left[ \begin{array}{cc}
a & b\\
c & d
\end{array} \right] \in SL(2, \rrr)$.

        We emphasize that the action of $SL(2, \rrr)$ on $\hhh$ is
{\em transitive}, i.e., for any two points $w_1, w_2 \in \hhh$ there
is an
element $g \in SL(2, \rrr)$ such that $w_2 = g\cdot w_1$. This can
easily
be proved. We also emphasize that there are subgroups of $SL(2, \rrr)$
for which the action is not transitive, among them the class of arithmetic
subgroups are to be mentioned. For example, the group $SL(2, \zzz)$ does
not act transitively on $\hhh$, and the set of orbits of the action of
$SL(2, \zzz)$ on $\hhh$, and similarly any arithmetic subgroup, is
infinite. We call the {\em arithmetic quotient} $\Gamma \backslash \hhh$
the set of orbits of the action of an arithmetic subgroup $\Gamma$ on
$\hhh$.

\begin{example} Take $\Gamma$ to be the {\em Hecke subgroup} $\Gamma_0(N)$
defined by
\[
\Gamma_0(N) = \left\{ \left[ \begin{array}{cc}
a & b\\
c & d
\end{array} \right] \in SL(2, \zzz)\;|\; c \equiv 0 ({\rm mod}\, N)\right
\},
\]
for a natural number $N$. This is a congruence subgroup and
$Y_0(N)=\Gamma_0(N)\backslash \hhh$ is an arithmetic
quotient. Such a quotient is not a compact subset, nor a bounded one, it is
however a subset with finite measure (volume) under the non-Euclidean measure
induced on the quotient from the group $SL(2, \rrr)$ which is a locally
compact group and induces the invariant volume element $\frac{dx \wedge
dy}{y^2}$, where $x, y$ are the real and the complex part of an element
$z\in \hhh$.
\end{example}

We now recall the basic ideas that turns an arithmetic quotient of the
form $\Gamma \backslash \hhh$ into an algebraic curve. Let $\Gamma \subset
SL(2, \qqq)$ be an arithmetic subgroup. The topological
boundary of $\hhh$ is $\rrr$ and a point $\infty$. For the rational
compactification of $\hhh$ we do not need to consider all the boundaries
$\rrr$ and $\{\infty\}$. In fact we need only to  add to $\hhh$ the cusps
of $\Gamma$ (a {\bf cusp} of $\Gamma$ is a rational number (an element of
$\qqq$) that is fixed under the action of an element $\gamma$ with the
property that $|tr(\gamma)| = 2$). Any two cusps $x_1, x_2$ such that
$\delta\cdot x_2 = x_1$ for an element $\delta \in \Gamma$  are called
{\bf equivalent}. Let $C(\Gamma)$ be the  set of inequivalent cusps of
$\Gamma$. Then $C(\Gamma)$ is finite. We add this set to $\hhh$ and form
the space $\hhh^\ast =\hhh \cup C(\Gamma)$. This space will be equipped
with certain topology such that a basis of the neighborhoods  of
the points of $\hhh^\ast$ is given by three type of open sets; if a point
in $\hhh^\ast$ is lying in $\hhh$ then its neighborhoods consists of usual
open discs in $\hhh$, if the point is $\infty$, i.e., the cusp $\infty$,
then its neighborhoods are the set of all points lying above the line
$Im(z) > \alpha$ for any real number $\alpha$, if the point is a cusp
different than $\infty$ which is a rational number, then the system of
neighborhoods of this point are the union of the cusp and the  interior of
a circle in $\hhh$ tangent to the cusp. Under the topology whose system of
open neighborhoods we just explained, $\hhh^\ast$ becomes a Hausdorff
non-locally compact space. The quotient space $\Gamma\backslash \hhh^\ast$
with the quotient topology is a compact Hausdorff space. We refer to this
compact quotient as the {\bf rational compactification} of
$\Gamma \backslash \hhh$. For a detailed discussion we refer the reader to
\cite{Shim1}.

        When the arithmetic group is a congruence subgroup of $SL(2, \zzz)$
the resulting algebraic curve is called a {\bf modular curve}.
For example, the rational compactification of $Y(N)=\Gamma(N) \backslash \hhh$
is denoted by $X(N)$ and the compactification of
$Y_0(N)=\Gamma_0(N) \backslash \hhh$ by $X_0(N)$.

\begin{example}
Let $N = 1$. Then $\Gamma =\Gamma(1) = SL(2, \zzz)$. In
this case $C(\Gamma) = \{\infty\}$, since all rational cusps are equivalent 
to the cusp $\infty$. So $\hhh^\ast = \hhh \cup \{\infty\}$, and
$\Gamma \backslash \hhh^\ast$ will be identified by
$\Gamma \backslash \hhh \cup \{\infty\}$. This may be seen as adding
$\infty$ to the fundamental domain $\fff_1 = \fff$ of $SL(2, \zzz)$, that
consists of all complex numbers in $z \in \hhh$ with
$|z| \geq1$ and $|Re(z)| \leq \frac{1}{2}$.
\end{example}

        The rational compactification of $\Gamma \backslash \hhh$ turns the
space $\Gamma \backslash \hhh^\ast$ into a compact Riemann surface (cf.
\cite{Shim1}) and so into an algebraic curve (cf. \cite{Nara}, or
\cite{SS}).

In general it is easiest to work with those arithmetic subgroups
which are torsion free and we shall assume from this point on that the
arithmetic subgroups we deal with have this property.
For example $\Gamma (N)$ and $\Gamma_0(N)$ for $N \geq 3$ are
torsion free.

\subsection {An adelic view of arithmetic quotients}

Consider the number field $\qqq$, the field of rational numbers. Let $\qqq_p$ be
the completion of $\qqq$ under the $p$-adic absolute value $|...|_p$, where
$|a/b|_p= p^{-n}$ whenever $a,b$ are integers and
$a/b= p^n\ddd\prod_{\ell\not= p\ {\rm prime}}\ell^{e_\ell}$, $n,e_\ell\in
\zzz$.
Recall that under the ordinary absolute value the completion of $\qqq$ is
$\rrr$. The ring of adeles of $\qqq$ is the locally compact commutative
ring $\aaa $ that is given by:
\[
\aaa  = \{ (x_\infty, x_2, \cdots) \in \rrr \times \ddd\prod_p \qqq_p\;|\;
\;\mbox{all but a finite number of}\;\; x_p \in \zzz_p\},
\]
where $\zzz_p$ is the ring of integers of $\qqq_p$ (as it is well known
$\zzz_p$ is a maximal compact open subring of $\qqq_p$). An element of
$\aaa $ is called an {\bf adele}. If $\aaa _f$ denotes the set of
adeles omitting the $\rrr$-component $x_\infty$, then $\aaa _f$ is
called the {\bf ring of finite adeles} and we can write
$\aaa  =3D \rrr \times \aaa _f$. Under the diagonal embedding $\qqq$ 
is a discrete subgroup of $\aaa $.

        We now consider the group $G = GL(2)$. For a choice of an
open compact subgroup $K_f \subset G(\aaa _f)$,
it is known that we can write the arithmetic
quotient (which was originally attached to an arithmetic subgroup of $\Gamma
\subset SL(2, \qqq)$) as the following quotient
\begin{equation}
Y(K_f) = G(\qqq) \backslash [\hhh \times (G(\aaa _f) /K_f)]=
\Gamma\backslash H,
\end{equation}
where
\begin{equation}
\Gamma = G(\qqq) \cap G(\rrr) K_f.
\label{eqn:Gamma}
\end{equation}
Thus our
arithmetic subgroup $\Gamma$ is completely determined by $K_f$. From now on
we assume $K_f$ has been chosen so that
$\Gamma$ is torsion free.

\begin{definition}
Let $G = GL(2)$. To $G$ is associated the Shimura
variety $Sh(G)$ as follows. Let $N \geq 3$ be a natural number. Let $\Gamma(N)$
be the congruence subgroup of level $N$ of $SL(2, \zzz)$, and $K = SO(2
, \rrr)$
the orthogonal group of $2 \times 2$ real matrices $A$ with determinant $1$
satisfying
$\phantom{.}^tA {A} = I_2$. Then
\[
Y(N)=\Gamma(N) \backslash \hhh \cong \Gamma(N) \backslash G(\rrr) / K.
\]
We call this the {\em modular space of level $N$}. Let
\[
K_f(N) = \{ g \in G(\ddd\prod_{p}\zzz_p)\;|\; g \equiv I_{2}({\rm mod}\, N)
\}
\]
be the {\bf open compact subgroup of $G(\aaa _f)$ of level $N$}. Then the
modular space of level $N$ can be written as:
\[
Y(N) \cong G(\qqq) \backslash G(\aaa )/K K_f(N) =
G(\qqq) \backslash [\hhh \times (G(\aaa _f)/K_f(N))].
\]
Thus
\[
X(K_f(N)) \cong Y(N).
\]
Taking the projective limit over $K_f(N)$ by letting $N$ gets large (which means
$K_f(N)$ gets small), we see that $\lim_N Y(N)
= G(\qqq)\backslash [\hhh \times G(\aaa _f)]$.
Then the (complex points of the) {\bf Shimura curve} $Sh(G)$ associated
to $G=SL(2)$ is defined by
\begin{equation}
Sh(G)(\ccc) = G(\qqq) \backslash [\hhh \times G(\aaa _f)].
\end{equation}
\end{definition}

Many mathematicians have addressed the natural questions
\begin{itemize}
\item
 What field are the curves $X(N)$, $X_0(N)$ defined over?
\item
 How can they be described  explicitly using
algebraic equations?
\end{itemize}

Regarding the first question, by the general theory of Shimura varieties we
know that for each reductive group $G$ defined over $\qqq$
satisfying the axioms of \S 2.1.1 in \cite{Del},
there is an algebraic number field $E=E_G$ over which a
Shimura variety $Sh(G)$ is defined \cite{Del}. In fact, the
Shimura curves $X(N)$
and $X_0(N)$ are regular schemes proper over $\zzz[1/N]$ (more precisely
over $Spec(\zzz[1/N])$)
\footnote{This result was essentially first proved by
Igusa \cite{Ig} (from the classical perspective).
See also \cite{TV}, Theorem 4.1.48,
\cite{Cas1} for an interesting discussion
of what happens at the ``bad primes'', and Deligne's paper in
the same volume as \cite{Cas1}.}.

Regarding the second question,
it is possible to find a {\bf modular polynomial} $H_N(x,y)$ of degree

\[
\mu(N)=N\ddd\prod_{p|N}(1+{\frac{1}{p}})
\]
for which $H_N(x,y)= 0$ describes (an affine patch of) $X_0(N)$.
Let
\[
G_k(q)=2\,\zeta (k)+2\, {\frac{( 2\,i\pi )^{k}}{(k-1)!}}
\ddd\sum _{n=1}^{\infty}\sigma_{k-1}(n)q^n,
\]
where $q = e^{2\pi i z}, z \in \hhh$, $\sigma_r(n)=\ddd\sum_{d|n}d^r$, and let
\[
\Delta(q)=60^3G_4(q)^3-27\cdot 140^2G_6(q)^2
=q\ddd\prod_{n=1}^\infty (1-q^n)^{24}.
\]
Define the {\bf j-invariant} by

\[
j(q)=1728\cdot 60^3G_4(q)^3/\Delta(q)
= q^{-1}+744+196884q+21493760q^2+864299970q^3+... \ .
\]
(More details on $\Delta$ and $j$ can be found
for example in \cite{Shim1}.)
The key property satisfied by $H_N$ is $H_N(j(q),j(q^N))= 0$.
It is interesting to note in passing
that when $N$ is such that the genus of $X_0(N)$ equals 0
(i.e., $N\in \{1, 3, 4, 5, 6, 7, 8, 9, 12, 13, 16, 18, 25\}$
\cite{Kn}) then this implies that $(x,y)= (j(q),j(q^N))$ parameterizes
$X_0(N)$. In general, comparing $q$-coefficients
allows one to compute $H_N$ for relatively small values of $N$.
(The MAGMA command {\tt ClassicalModularPolynomial}\footnote{See 
also {\tt CanonicalModularPolynomial} and
{\tt AtkinModularPolynomial},} computes this expression.
However, even for $N= 11$, some of the coefficients
can involve one hundred digits or more.
The cases $N= 2,3$ are given in Elkies \cite{E}, for example.
The paper by P. Cohen \cite{Co} determines the asymptotic
size of the largest coefficient of $H_N$ (normalized to have
leading coefficient equal to $1$). She shows that the
largest coefficient grows like $N^{c\mu(N)}$,
where $c>0$ is a constant.
More practical equations for (some of) the $X_0(N)$ are given in
T. Hibino and
N. Murabayashi \cite{HM},
M. Shimura \cite{Shim2}, J. Rovira \cite{R}, G. Frey and M. M\"uller
\cite{FM}, Birch \cite{B}, and the table in \S 2.5 below.

For deeper study of Shimura varieties and the theory of
canonical models we refer the reader to \cite{Del},
 \cite{Lan2}, and \cite{Shim1}.

\subsection{Hecke operators and arithmetic on $X_0(N)$}

In this section we recall some well-known though relatively deep
results on $X_0(N)(\fff_p)$, where $p$ is a prime not dividing $N$.
These shall be used in the discussion of the Tsfasman, Vladut,
Zink, and Ihara result later.

First, some notation: let $S_2(\Gamma_0(N))$ denote the space of
holomorphic automorphic forms of weight 2 on $\Gamma_0(N)\backslash H$.
Let $T_p:S_2(\Gamma_0(N))\rightarrow S_2(\Gamma_0(N))$ denote the
{\bf Hecke operator} defined by

\[
T_pf(z)=
f(pz)+\ddd\sum_{i= 0}^{p-1}
f({\frac{z+i}{p}}),\ \ \ \ z\in H.
\]
Define $T_{p^k}$ inductively by
\[
T_{p^k}= T_{p^{k-1}}T_{p}-pT_{p^{k-2}},\ \ \ T_1=1,
\]
and define the modified Hecke operators $U_{p^k}$ by
\[
U_{p^k}=T_{p^{k}}-pT_{p^{k-2}},\ \ \ U_p= T_p,
\]
for $k\geq 2$. The Hecke operators may be extended to the positive integers
by demanding that they be multiplicative.
\begin{theorem}
\label{thrm:Eichler-Shimura}
(``Congruence relation'' of Eichler-Shimura
 \cite{M}, \S 5.6.7, or \cite{St}) Let $q=p^k$, $k>0$ an integer.
If $p$ is a prime not dividing $N$ then

\[
Tr(T_p)=p+1-|X_0(N)(\fff_p)|.
\]
More generally,

\[
Tr(T_q-pT_{q/p^2})=q+1-|X_0(N)(\fff_q)|.
\]
\end{theorem}

\begin{example}
One may try to compute the trace of the Hecke operators
$T_p$ acting on the space of holomorphic
cusp forms of weight $2$, $S_2(\Gamma_0(N))$, by using either
the Eichler-Shimura congruence relation, which we give
below (see Theorem \ref{thrm:Eichler-Shimura}),
or by using some easier but ad hoc ideas going back to
Hecke which work in special cases. One simple idea is to
note that $S_2(\Gamma_0(N))$ is spanned by
simultaneous eigenforms of the Hecke operators
(see for example, Proposition 51 in chapter III of \cite{Ko}).
In this case, it is known that
the Fourier coefficient $a_p$, $p$
prime not dividing $N$, of a normalized
(to have leading coefficient $a_1=1$)
eigenform is the eigenvalue of $T_p$ (see
for example, Proposition 40 in chapter III of \cite{Ko}).
If $S_2(\Gamma_0(N))$ is one-dimensional then any
element in that space $f(z)$ is such an eigenform.

The modular curve $X_0(11)$ is of genus 1, so there is
(up to a non-zero constant factor) only one holomorphic
cusp form of weight $2$ in $S_2(\Gamma_0(11))$
(see Theorem \ref{thrm:HZformula} below).
There is a well-known construction of this form
(see \cite{O2} or \cite{Gel}, Example 5.1),
which we recall below.
As we noted above, the $p$-th coefficient $a_p$
($p$ a prime distinct from $11$)
of its Fourier expansion is known to satisfy
$a_p=Tr(T_p)$. These will be computed using MAGMA.

Let $q=e^{2\pi i z}$, $z\in \hhh$, and
consider {\bf Dedekind's $\eta$-function}:
\index{Dedekind $\eta$-function}
\[
\eta(z)=e^{2\pi i z/24}\ddd\prod_{n=1}^\infty (1-q^n).
\]
Then
\[
f(z)=\eta(z)^2\eta(11z)^2
q\ddd\prod_{n=1}^\infty (1-q^n)^{2}(1-q^{11n})^{2},
\]
is an element of $S_2(\Gamma_0(11))$
\footnote{In fact, if we write $f(z)=\ddd\sum_{n=1}^\infty a_nq^n$ then
\[
\zeta_E(s)=(1-p^{-s})^{-1}
\ddd\prod_{p\not= 11}(1-a_pp^{-s}+p^{1-2s})^{-1},
\]
is the global Hasse-Weil zeta function of the elliptic curve
$E$ of conductor $11$
with Weierstass model $y^2+y=x^3-x^2$ \cite{Gel}, page 252.}.
One can compute the $q$-expansion of this form using MAGMA's
{\tt ModularForms(Gamma0(11),2)} command:

\[
f(z)=q-2\,{q}^{2}-{q}^{3}+2\,{q}^{4}+{q}^{5}+2\,{q}^{6}-2\,{q}^{7}...
\]
For example, the above expansion tells us that $Tr(T_3)= Tr(U_3)=-1$.
The curve $X_0(11)$ is of genus 1 and is
isogenous to the elliptic curve $E$
with Weierstrass model $y^2+y=x^3-x^2$. Over the field with $p=3$ elements,
there are $|X_0(11)(\fff_3)|=p+1-Tr(T_p)=5$ points in $E(\fff_3)$,
including $\infty$:
\[
E(\fff_3)=\{[0, 0], [0, 2], [1, 0], [1, 2],\infty\}.
\]
(For this, one uses the commands
\verb+F:=GF(3); P<x>:=PolynomialRing(F);+
\verb+f:=x^3-x^2;h:=1; C:=HyperellipticCurve(f, h);Places(C,1);+.)

For a representation-theoretic discussion of this
example, see \cite{Gel}, \S 14.

For an example of an explicit element of $S_2(\Gamma_0(32))$,
see Koblitz \cite{Ko}, chapter II \S 5 and (3.40) in chapter III.
For a remarkable theorem which illustrates how far this
$\eta$-function construction can be extended, see Morris' theorem
in \S 2.2 of \cite{R}.
\end{example}

To estimate $a_{p^{k}}$, one may appeal to an explicit
expression for $Tr(T_{p^{k}})$ known as the
``Eichler-Selberg trace formula'', which we discuss next.

\subsection{Eichler-Selberg trace formula}

In this subsection, we recall the version of the trace formula for
the Hecke operators due to Duflo-Labesse \cite{DL}, \S 6.

Let $k$ be an even positive integer and
let $\Gamma$ be a congruence subgroup as in (\ref{eqn:Gamma}).
Let $S$ denote a complete set of representatives of
$G(\qqq)$-conjugacy classes of $\rrr$-elliptic elements in $\Gamma$
($\rrr$-elliptic elements are those that are conjugate to an element of
$SO(2, \rrr)$, the orthogonal group). For $\gamma\in S$, let $w(\gamma)$
denote the cardinality of the centralizer of $\gamma$ in $\Gamma$. If
$r(\theta)
=\left(
\begin{array}{cc}
\cos(\theta)&\sin(\theta)\\
-\sin(\theta)&\cos(\theta)
\end{array}
\right)
$ then let $\theta_\gamma \in (0,2\pi)$ denote the element for which
$\gamma= r(\theta_\gamma)$. Let $\tau_m$ denote the
image in $G(\aaa_f)$ of the set of matrices
in $GL(2,\aaa_f)$ having coefficients in
$\hat{\zzz}= \ddd\prod_{p<\infty}\zzz_p$ and determinant in
$m\hat{\zzz}$.
Consider the subspace $S_k(\Gamma)\subset L^2(\Gamma\backslash H)$
formed by functions satisfying
\begin{itemize}
\item
$f(\gamma z)= (cz+d)^{k}f(z)$, for all
$\gamma
=\left(
\begin{array}{cc}
a&b\\
c&d
\end{array}
\right)\in \Gamma$, $x\in H$,
\item
$f$ is a holomorphic cusp form.
\end{itemize}
This is the space of holomorphic cusp forms of weight $k$ on $\hhh$.

Let
\[
\epsilon(\sqrt{m})=
\left\{
\begin{array}{cc}
1,&  m\ {\rm is\ a\ square},\\
0, & {\rm otherwise}.
\end{array}
\right.
\]
and let
\[
\delta_{i,j}=
\left\{
\begin{array}{cc}
1,&  i= j,\\
0, & {\rm otherwise}.
\end{array}
\right. .
\]

\begin{theorem}
(``Eichler-Selberg trace formula'')
\index{Eichler-Selberg trace formula}
\label{thrm:Eichler-Selberg}
Let $k>0$ be an even integer and $m>0$ an integer.
The trace of $T_m$ acting on $S_k(\Gamma)$ is given by
\[
\begin{array}{c}
Tr(T_m)= \delta_{2,k}\ddd\sum_{d|m}b +
\epsilon(\sqrt{m})({\frac{k-1}{12}}m^{(k-2)/2}-{\frac{1}{2}}m^{(k-1)/2})\\
-\ddd\sum_{\gamma\in S\cap \tau_m}w(\gamma)^{-1}m^{(k-2)/2}
{\frac{\sin((k-1)\theta_\gamma)}{\sin(\theta_\gamma)}}
-\ddd\sum_{d|m,\ d^2<m}b^{k-1}.
\end{array}
\]
\end{theorem}

\begin{remark}
Let $k= 2$, $m = p^2$, $\Gamma= \Gamma_0(N)$
and $N\rightarrow \infty$ in the above formula.
It is possible to show that the Eichler-Selberg trace formula
implies
\begin{equation}
Tr(T_{p^2})= g(X_0(N))+O(1),
\label{eqn:estimate}
\end{equation}
as $N\rightarrow \infty$. The proof of this
estimate (see \cite{M}, chapter 5, or \cite{LvdG},
\S V.4)
uses the explicit formula given below for $g(X_0(N))=
{\rm dim}(S_2(\Gamma_0(N))$, which we shall also make use of
later.

\begin{theorem} (``Hurwitz-Zeuthen formula'' \cite{Shim1})
\index{Hurwitz-Zeuthen formula}
\label{thrm:HZformula}
\footnote{The genus formula for $X_0(N)$ given in
\cite{Shim1} and \cite{Kn} both apparently contain a
(typographical?) error. The problem is in the
$\mu_2$ term, which should contain
a Legendre symbol $({\frac{-4}{n}})$ instead of
$({\frac{-1}{n}})$. See for example \cite{Ei} for a
correct generalization.}
The genus of $X_0(N)$ is given by
\[
g(X_0(N))={\rm dim}(S_2(\Gamma_0(N))=
1+{\frac{1}{12}}\mu(N)-{\frac{1}{4}}\mu_2(N)-
{\frac{1}{3}}\mu_3(N)-\mu_\infty(N),
\]
where
\[
\mu(N)=[SL(2, \zzz)/\Gamma_0(N)]=N\ddd\prod_{p|N}(1+{\frac{1}{p}}),
\]
\[
\mu_2(N)=
\left\{
\begin{array}{cc}
\ddd\prod_{p|N\ {\rm prime}}(1+({\frac{-4}{p}})),& gcd(4,N)=1,\\
0,&4|N,
\end{array}
\right.
\]
\[
\mu_3(N)=
\left\{
\begin{array}{cc}
\ddd\prod_{p|N\ {\rm prime}}(1+({\frac{-3}{p}})),& gcd(2,N)=1\ {\rm and}\
gcd(9,N)\not= 9,\\
0,&2|N\ {\rm or}\ 9|N,
\end{array}
\right.
\]
and
\[
\mu_\infty(N)=\ddd\sum_{d|N}\phi(gcd(d,N/d)),
\]
where $\phi$ is Euler's totient function and $({\frac{\cdot}{p}})$ is
Legendre's symbol.
\end{theorem}

The estimate (\ref{eqn:estimate}) and the
Eichler-Shimura congruence relation imply
\[
\begin{array}{c}
|X_0(N)(\fff_{p^2})|= p^2+1-Tr(T_{p^2}-pI)= p^2+1-Tr(T_{p^2})+
p\cdot {\rm dim}(S_2(\Gamma_0(N))\\
= p^2+1-(g(X_0(N))+O(1)))+p\cdot g(X_0(N))\\
= (p-1)g(X_0(N))+O(1),
\end{array}
\]
as $N\rightarrow \infty$.
\end{remark}

\subsection{The curves $X_0(N)$ of genus 1}

It is known (see for example \cite{Kn})
that a modular curve of level $N$, $X_0(N)$, is of genus 1
if and only if
\[
N\in \{11,14,15,17,19,20,21,24,32,36,49\}.
\]
In these cases, $X_0(N)$ is birational to an elliptic curve
$E$ having Weierstrass model of the form
\[
y^{2}+ a_1xy+a_3y= {x}^{3}+a_2x^2 +a_4 x +a_6,
\]
with $a_1,a_2,a_3,a_4,a_6$.
If $E$ is of above form then the {\bf discriminant} is given by
\index{discriminant}
\[
\Delta=
-b_2^2b_8-8b_4^3-27b_6^2+9b_2b_4b_6,
\]
where
\[
b_2= a_1^2+4a_2,\ \ \ \ \
b_4=2a_4+a_1a_3,\ \ \ \ \
b_6=a_3^2+4a_6,
\]
\[
b_8=a_1^2a_6+2a_2a_6-a_1a_3a_4+a_2a_3^2-a_4^2.
\]
The conductor
\footnote{The conductor is defined in Ogg \cite{O1}, but see also
\cite{Gel}, \S I.2, or \cite{Kn}, P. 390.}
$N$ of $E$ and its discriminant $\Delta$ have the same prime factors.
Furthermore, $N|\Delta$ (\cite{Kn}, \cite{Gel}).

Some examples, which we shall use later,
are collected in the following table.

\begin{table}
\label{table:models}
\begin{center}
\begin{tabular}{|c|c|c|c|}\hline
level & discriminant & Weierstrass model & reference \\ \hline
11 & -11 & ${y}^{2}+y={x}^{3}-{x}^{2}$ & \cite{BK}, table 1, p. 82 \\
14 &-28 & ${y}^{2}+xy-y={x}^{3}$ &p. 391, table 12.1 of \cite{Kn} \\
15 &15 & ${y}^{2}+7\,xy+2\,y={x}^{3}+4\,{x}^{2}+x$ & p. 65, table 3.2
of \cite{Kn}\\
17 &17 & ${y}^{2}+3\,xy={x}^{3}+x$ & p. 65, table 3.2 of \cite{Kn} \\
19  & -19 & ${y}^{2}+y={x}^{3}+{x}^{2}+x$ &\cite{BK}, table 1, p. 82 \\
20 &80 & ${y}^{2}={x}^{3}+{x}^{2}-x$ & p. 391, table 12.1 of \cite{Kn}\\
21 & -63 & ${y}^{2}+xy={x}^{3}+x$ &p. 391, table 12.1 of \cite{Kn} \\
24 &-48 & ${y}^{2}={x}^{3}-{x}^{2}+x$ & p. 391, table 12.1 of
\cite{Kn}\\
27  &-27 & ${y}^{2}+y={x}^{3}$ & p. 391, table 12.1 of \cite{Kn}\\
32  &64 & ${y}^{2}={x}^{3}-x$ & p. 391, table 12.1 of \cite{Kn}\\
36  & & (see below) &\S 4.3 in \cite{R} \\
49  &  &(see below) & \S 4.3 in \cite{R}\\ \hline
\end{tabular}
\caption{Models of genus $1$ modular curves}
\end{center}
\end{table}

When $N=36$, \S 4.3 in Rovira \cite{R} gives ${y}^{2}={x}^4-4x^3-6x^2-4x+1$,
which is a hyperelliptic equation but not in Weierstrass form.
To put it in Weierstrass form, we use \cite{MAGMA}\footnote{More 
precisely, we use the
{\tt ReducedMinimalWeierstrassModel} command over the
field $\qqq$.}.
This produces the cubic equation 
$y^2 + (x^2 + 1)y = x^3 - 2x^2 + x$, provided
$p\not= 2$. This has conductor $\Delta=-1769472$.
When $N=49$, \S 4.3 in Rovira \cite{R} gives $y^2=x^4-2x^3-9x^2+10x-3$,
which is a hyperelliptic equation but not in Weierstrass form.
As before, MAGMA produces a cubic equation in which the coefficient of
$x^3$ is not one,
$y^2 + (-x^2 - x - 1)y = -x^3 - 3x^2 + 2x - 1$.
The change-of-variable $x\longmapsto -x$
produced the Weierstrass form
$y^2 + (-x^2 - x - 1)y = x^3 - 3x^2 - 2x - 1$.
This has conductor $\Delta=-1404928$.

\section{Codes}

        To have an idea of how the points of a curve over finite fields are
used in the coding theory we first recall the definition of a code.

        Let $A$ be a finite set, which we regard as an
alphabet. Let $A^n $ be the
$n$-fold Cartesian product of $A$ by itself. In $A^n$ we define the
{\bf Hamming metric} $d(x, y)$ by:
\index{Hamming!metric}\index{Hamming!distance}
\[
d(x, y) = d((x_1, \cdots, x_n), (y_1, \cdots, y_n)) : = |\{i\;|\;
x_i \not= y_i\}|.
\]
We now assume that $A^n$ is equipped with the Hamming metric. Then by
definition a subset $C \subseteq A^n$ is called an
{\bf $|A|$-ary code}. An important case
\index{Code!(definition of error-correcting ...)}
arises when we let $A$ to be a finite field.
Suppose that $q = p^m$ and $\fff_q$ is a finite
field with $q$ elements. In this case we may put
$A = \fff_q$ and $y = (0, \cdots,0)$.
Then the {\bf weight} of $x$ is the Hamming length
\index{weight of a codeword}
$||x|| = d(x, 0) = |\{i\;|\; x_i \not=0\}|$. In
particular a subset $C$ of $\fff^n_q$ is a code,
and to it is associated two basic
parameters: $k = \log_q|C|$, the {\bf number of information bits}
\index{information bits of a codeword}
and $d = min\{||x - y||\;|\; x, y \in C, y\not=0\}$ the {\bf minimum 
distance}.
\index{minimum distance of a code}
(A code with minimum distance $d$ can correct $[{\frac{d-1}{2}}]$
errors.)
Let
\[
R = R(C) = \frac{k}{n},
\]
which measures the information rate of the code, and
\[
\delta = \delta(C) = \frac{d}{n},
\]
which measures the error correcting ability of the code.

\subsection{Basics on linear codes}

If the code $C\subset \fff^n_q$ is a vector space over $\fff_q$ then
we call $C$ a {\bf linear code}. The {\bf parameters}
\index{linear code}\index{parameters of a code}
of a linear code $C$ are
\begin{itemize}
\item the {\bf length} $n$,
\index{length of a code}
\item
the {\bf dimension} $k={\rm dim}_{\fff_q}(C)$,
\index{dimension of a linear code}
\item
the minimum distance $d$.
\end{itemize}
Such a code is called an {\bf $(n,k,d)$-code}.
Let $\Sigma_q$ denote the set of all $(\delta,R)\in [0,1]^2$ such that
there exists a sequence $C_i$, $i=1,2,...$, of $(n_i,k_i,d_i)$-codes
for which $\lim_{i\rightarrow \infty} \delta_i=\delta$
and $\lim_{i\rightarrow \infty} R_i=R$.

The following theorem describes
information-theoretical limits on how ``good''
a linear code can be.

\begin{theorem} (Manin \cite{SS}, chapter 1)
\index{Manin's theorem on code rates}
There exists a continuous decreasing function
\[
\alpha_q:[0,1]\rightarrow [0,1],
\]
such that
\begin{itemize}
\item
$\alpha_q$ is strictly decreasing on $[0,{\frac{q-1}{q}}]$,
\item
$\alpha_q(0)=1$,
\item
if ${\frac{q-1}{q}}\leq x\leq 1$ then $\alpha_q(x)=0$,
\item
$\Sigma_q=\{(\delta,R)\in [0,1]^2\ |\ 0\leq R\leq \alpha_q(\delta)\}$.

\end{itemize}
\end{theorem}

Not a single value of $\alpha_q(x)$ is known for $0<x<{\frac{q-1}{q}}$!
It is not known whether or not the maximum
value of the bound, $R= \alpha_q(\delta)$
is attained by a sequence of linear codes.
It is not known whether or not
$\alpha_q(x)$ is differentiable for $0<x<{\frac{q-1}{q}}$,
nor is it known if $\alpha_q(x)$ is convex on
$0<x<{\frac{q-1}{q}}$.
However, the following estimate is known.

\begin{theorem} (Gilbert-Varshamov \cite{MS}, \cite{SS} chapter 1)
\index{Gilbert-Varshamov!bound}
\label{thrm:GV}
We have
\[
\alpha_q(x)\geq 1- x\log_q(q-1)-x\log_q(x)-(1-x)\log_q(1-x).
\]
In other words, for each fixed $\epsilon >0$,
there exists an $(n,k,d)$-code $C$ (which may depend on $\epsilon$)
with
\[
R(C)+\delta(C)\geq
1- \delta(C)\log_q({\frac{q-1}{q}})-\delta(C)\log_q(\delta(C))-
(1-\delta(C))\log_q(1-\delta(C))-\epsilon.
\]
\end{theorem}

The curve $(\delta, 1- \delta\log_q({\frac{q-1}{q}})-\delta\log_q(\delta)-
(1-\delta)\log_q(1-\delta)))$ is called the
{\bf Gilbert-Varshamov curve}.
This theorem says nothing about constructing codes satisfying this
property! Nor was it known, until the work of
Tsfasman, Valdut, Zink and Ihara, how to do so.

\subsection{Some basics on AG codes}

We begin with Goppa's basic idea boiled down to its most basic
form. Let $R$ denote a commutative ring with unit and
let $m_1,m_2,...,m_n$ denote a
finite number of maximal ideals such that
for each $1\leq i\leq n$, we have $R/m_i\cong \fff_q$.
Define $\gamma:R\rightarrow \fff_q^n$ by
\[
\gamma(x)=(x+m_1,x+m_2,...,x+m_n),\ \ \ \ \ x\in R.
\]
Of course, in this level of generality, one cannot say much
about this map. However, when $R$ is associated to
the coordinate functions of a
curve defined over $\fff_q$ then one can often use the machinery
of algebraic geometry to obtain good
estimates on the parameters $(n,k,d)$ of the
code associated to $\gamma$.

Let $V$ be an irreducible smooth projective
algebraic variety defined over the finite field $\fff_q$.
Let $\fff_q(V)$ denote the field of rational functions
on $V$.
Let ${\cal P}(V)$ denote the set of prime divisors of $V$,
which we may identify with the closed
irreducible subvarieties of
$V(\overline{\fff_p})$
of codimension $1$.
For each $P\in {\cal P}(V)$, there is a valuation
map $ord_P:\fff_q(V)\rightarrow \zzz$
(see Hartshorne \cite{Ha}, \S II.6, page 130).
Let ${\cal D}(V)$ denote the group of divisors of $V$,
the free abelian group generated by ${\cal P}(V)$.
\index{divisor group (of a variety)}

If $A=\ddd\sum_P a_P P, B=\ddd\sum_P b_P P\in {\cal D}(V)$
are divisors then we say $A\leq B$ if and only if
$a_P\leq b_P$ for all $P\in {\cal P}(V)$. If
$f\in \fff_q(V)$ is a non-zero function then let
\[
div(f)=\ddd\sum_{P\in {\cal P}(V)} ord_P(f)P,
\]
where $ord_P(f)$ is the order of the zero
(pole) at $P$ (as above).
This is well-defined (since the above sum is
finite by Lemma 6.1 in \cite{Ha}, \S II.6, page 131).
For $B\in {\cal D}(V)$, define ${\cal L}(B)=H^0(V,{\cal O}_B)$ 
to be the Riemann-Roch space

\[
{\cal L}(B)=\{0\}\cup \{f\in \fff_q(V)\ |\ f\not= 0, div(f)\geq -B\}.
\]

Pick $n$ different points $P_1,P_2,...,P_n$
in $V(\fff_q)$, let $D=P_1+...+P_n$, 
and choose a divisor
$G=\ddd\sum_{P\in {\cal P}(V)} a_P P\in {\cal D}(V)$
disjoint from these points (i.e., no $P_i$ is
a point on the codimension one subvariety $P$ in $G$).
It is {\it not} necessary for $G$ to be rational.
The {\bf Goppa code} or {\bf AG code}
associated to $(V(\fff_q),D,G)$ is the
linear code $C=C(G,D,V)$ defined to
be the subspace of $\fff_q^n$ which is the image of
the map
\begin{equation}
\gamma :{\cal L}(G)\rightarrow \fff_q^n,
\label{eqn:gamma}
\end{equation}
defined by $\gamma(f)=(f(P_1),...,f(P_n))$. (In 
the case of curves, this code
is called the {\bf dual Goppa code}
or {\bf Goppa function code} in \cite{P}. Goppa gave another
geometric construction of codes using differentials
for which we refer the reader to \cite{P} or
\cite{TV}. In other parts of the literature, the term
``Goppa code'' refers to an earlier construction of Goppa 
using rational functions.)

To specify an AG code, one must

\begin{itemize}
\item
choose a smooth variety $V$ over $\fff_q$,
\item
pick rational points $P_1,P_2,...,P_n$ of $V$,
\item
choose a divisor $G$ disjoint from the $P_i$'s,
\item
determine a basis for ${\cal L}(G)$,
\item
compute the matrix for $\gamma$ wth respect to this basis.
\end{itemize}

\subsection{Some estimates on AG codes}

Let $g $ be the genus of a curve $V=X$ and let
$C=C(G,D,X)$ denote the Goppa code as constructed above.
If $G$ has parameters $[n,k,d]$ and if we then the following
lemma is a consequence of the Riemann-Roch theorem.

\begin{lemma} Assume $C$ is as above and $G$ satisfies 
$2g-2<deg(G)<n$.
Then $k=dim(C)=deg(G)-g+1$ and $d\geq n-deg(G)$.
\end{lemma}

Consequently, $k+d\geq n-g+1$. Because of Singleton's
inequality\footnote{It is known that $n\geq d+k-1$ for any linear
$(n,k,d)$-code (this is the {\bf Singleton inequality}),
\index{Singleton inequality}
with equality if and only if the code is a so-called
{\bf MDS code} (MDS=minimum distance separable).}, we have
\index{MDS!code}
\begin{itemize}
\item
if $g=0$ then $G$ is an MDS code,
\item
if $g=1$ then $n\leq k+d\leq n+1$.
\end{itemize}
The previous lemma also implies the following lower bound.

\begin{proposition} (\cite{SS} \S 3.1, or \cite{TV})
\label{prop:7}
With $C$ as in the previous lemma, we 
have $\delta +R= {\frac{d}{n}} + {\frac{k}{n}} \geq 1 - \frac{g-1}{n}$.
\end{proposition}

Theorem \ref{thrm:HZformula} above is an explicit formula for the
genus of the modular curve $X_0(N)$.
It may be instructive to
plug this formula into the estimate in Proposition \ref{prop:7}
to see what we get. The formula for the
genus $g_N$ of $X_0(N)$
is relatively complicated, but simplifies greatly
when $N$ is a prime number which is congruent to
$1$ modulo $12$, say $N=1+12m$, in which case
$g_N=m-1$. For example, $g_{13}=0$. In particular,
we have the following

\begin{corollary} Let $X=X_0(N)$, where
$N$ is a prime number which is congruent to
$1$ modulo $12$ and which has the property
that $X$ is smooth over $\fff_q$.
Then the parameters $[n,k,d]$ of a
Goppa code associated to $X$ must satisfy
\[
{\frac{d}{n}} + {\frac{k}{n}} \geq 1 - \frac{{\frac{N-1}{12}}-2}{n}.
\]
\end{corollary}

Based on the above Proposition, if one considers a family of curves
$X_i$ with increasing genus $g_i$ such that
\begin{equation}
\lim_{i\rightarrow \infty}{\frac{|X_i(\fff_q)|}{g_i}} = \alpha
\label{eqn:4}
\end{equation}
one can construct a family of codes $C_i$ with
$\delta(C_i) + R(C_i) \geq 1 - \frac{1}{\alpha}$.
It is known that $\alpha \leq \sqrt{q}-1$
(this is the so-called {\bf Drinfeld-Vladut bound},
\index{Drinfeld-Vladut bound}
\cite{TV}, Theorem 2.3.22).

The following result says that the Drinfeld-Vladut bound can
be attained in case $q=p^2$.

\begin{theorem} (Tsfasman, Valdut, Zink \cite{TV}, Theorem 4.1.52)
Let $g_N$ denote the genus of $X_0(N)$.
If $N$ runs over a set of primes different
than $p$ then the quotients $g_N/|X_0(N)(\fff_{p^2})|$
associated to the modular curves $X_0(N)$
tend to the limit $\frac{1}{p-1}$.
\end{theorem}

More generally, if
$q = p^{2 k}$, then there is a family of Drinfeld
curves $X_i$ over $\fff_q$
yielding $\alpha = \sqrt{q}-1$
(\cite{TV}, Theorem 4.2.38,
discovered independently by Ihara \cite{I} at about the same time). In other
words, the Drinfeld-Vladut bound is attained in case
$q = p^{2 k}$.

As a corollary to the above theorem,
if $p \geq 7$ then there exists a sequence of Goppa codes
$G_N$ over $\fff_{p^2}$
associated to a sequence of modular curves $X_0(N)$
for which $(R(G_N),\delta(G_N))$
eventually (for suitable large $N$)
lies above the Gilbert-Varshamov bound
in Theorem \ref{thrm:GV}. This follows from comparing the
Gilbert-Varshamov curve
\[
(\delta, 1- \delta\log_q({\frac{q-1}{q}})-\delta\log_q(\delta)-
(1-\delta)\log_q(1-\delta)))
\]
with the curve $(\delta,\frac{1}{\sqrt{q}-1})$, $q=p^2$.

\section{Examples}

Let $C$ be an elliptic curve. This is a projective
curve for which $C(\fff_q)$
has the structure of an algebraic
group. Let $P_0\in C(\fff_q)$ denote the identity.
Let $P_1,P_2,...,P_n$ denote {\it all}
the other elements of $C(\fff_q)$
and let $A=aP_0$, where $0<a<n$ is an integer.

\begin{example} Let $C$ denote the elliptic curve of conductor
$32$ (and birational to $X_0(32)$) with Weierstrass
form $y^2=x^3-x$. If $p$ is a prime
satisfying $p\equiv 3({\rm mod}\, 4)$ then
\[
|C(\fff_p)|=p+1
\]
(Theorem 5, \S 18.4 in Ireland and Rosen \cite{IR}).
Let $C(\fff_p)=\{P_0,P_1,P_2,...,P_n\}$,
where $P_0$ is the identity, and if $A=kP_0$, for some
$k>0$. The parameters of the corresponding code
$G=G(A,P,C)$ satisfy $n=p, \ \ d+k\geq n$,
since $g=1$, by the above Proposition.
As we observed above, an AG code constructed from an elliptic curve
satisfies either $d+k-1=n$ (i.e., is MDS) or else $d+k=n$.
The result of Shokrollahi below implies that
if, in addition, $p>3$ or $k>2$ then $G$ is not MDS and
\[
n=p,\ \ \ \ d+k=p.
\]
\end{example}

The following result is an immediate corollary of
the results in \cite{Sh}, see also \S 5.2.2 in \cite{TV}.

\begin{theorem} (Shokrollahi)
Let $C, P_0, P_1, ..., P_n,
D, A$, be as above.
\begin{itemize}
\item
If $a=2$ and $C(\fff_q)\cong C_2\times C_2$
(where $C_n$ denotes the cyclic group of order
$n$) then the code $G(A,D)$ is a
$[n,k,d]$-code ($n$ is the length, $k$ is the
dimension, and $d$ is the minimum distance)
with
\[
d=n-k+1,\ \ \ {\rm and}\ \ \  k=a.
\]
\item
Assume $gcd(n,a!)=1$.
If $a\not= 2$ or $C(\fff_q)$ is {\it not}
isomorphic to the Klein four group
$C_2\times C_2$ then $G(A,D)$ is a
$[n,k,d]$-code ($n$ is the length, $k$ is the
dimension, and $d$ is the minimum distance)
with
\[
k=a
\]
and weight enumerator polynomial (see for
example \cite{MS} for the definition)
\[
W_G(x)=x^n+\ddd\sum_{i=0}^{a-1}
\left(
\begin{array}{c}
n \\
i
\end{array}
\right)
(q^{a-i}-1)(x-1)^i +B_a(x-1)^a,
\]
where $B_a$ is given in \cite{Sh} and \S 3,2,2 in \cite{TV}.
\end{itemize}
\end{theorem}

\subsection{Weight enumerators of some elliptic codes}

In the case where $E$ is given by the level 19,
discriminant $-19$ Weierstrass model
\[
y^2+y=x^3+x^2+x,
\]
and $p=13$, we have
\[
E(\fff_p)=
\begin{array}{cc}
&\{\infty, [0, 0], [0, 12], [1, 6], [3, 0], [3, 12], [4, 2],
[4, 10], [5, 3], [5, 9],  \\
& [8, 3], [8, 9],
[9, 0], [9, 12], [11, 4], [11, 8], [12, 3], [12, 9]\}
\end{array}.
\]
Write $E(\fff_p)=
\{P_0,P_1,...,P_{17}\}$, where $P_0$ denotes the identity element
of the group law for $E$, let
$A=kP_0$, and let $D=P_1+...+P_{17}$.
The hypotheses of the above theorem are satisfied when
we take $n=17$ and $2\leq k=a<17$. The above construction
associates to this data a
Goppa code $G=G(A,D,E)$ which is a $7$-error correcting code of length
$n=17$ over $\fff_{13}$.
Some of the weight enumerator polynomials $W_G$ and the number
of errors these codes $G$ can correct are given in the following
table.

\begin{table}
\label{table:weights}
\begin{center}
\begin{tabular}{|c|c|c|} \hline
 &  & number of errors \\
$a=k$ & weight enumerator $W_G$ & $G$ corrects \\ \hline
$2$ & $x^{17}+96x^2+12x+60$ & $7$ \\
$3$ & $x^{17}+456x^3+264x^2+960x+516$ & $6$ \\
$4$ & $x^{17}+1608x^4+1728x^3+8016x^2+9684x+7524$ & $6$ \\
$5$ & $x^{17}+4104x^5+8040x^4+...+94644$ & $5$ \\
$6$ & $x^{17}+8232x^6+ 24864x^5+...+1239540$ & $5$ \\
$7$ & $x^{17}+12984x^7+57624x^6+...+16090116$ & $4$ \\
$8$ & $x^{17}+16272x^8+103200x^7+...+209219292$ & $4$ \\
$9$ & $x^{17}+16176x^9+146136x^8+...+2719777524$ & $3$ \\
$10$ & $x^{17}+12912x^{10}+162600x^9+...+35357193732$ & $3$ \\
\hline
\end{tabular}
\end{center}
\caption{weight enumerator polynomials of some elliptic codes}
\end{table}
These were computed using MAGMA commands such as the following.

\begin{verbatim}
F:=GF(13);
P<x>:=PolynomialRing(F);
f:=x^3+x^2+x;h:=1;
C:=HyperellipticCurve(f, h);
Places(C,1);
Div := DivisorGroup(C);
Pls:=Places(C,1);
S:=[Pls[i] : i in [2..#Pls]];
m:=2;
D := m*(Div!Pls[1]);
AGC := AlgebraicGeometricCode(S, D);
Length(AGC);
Dimension(AGC);
MinimumDistance(AGC);
WeightEnumerator(AGC);
\end{verbatim}

The number of codewords of minimum weight $n-k$ is the
coefficient of the second highest term in $W_G(x)$.
For example, when $k=3$ the
number of codewords of minimum weight $n-k=14$ is $384$.

A smaller example using the same elliptic curve
$E$ as above: taking $p=3$, we find
that
\[
E(\fff_p)=
\{
[0, 0], [0, 2], [1, 0], [1, 2], [2, 1],
\infty\}.
\]
The hypotheses of the above theorem are satisfied when
we take $n=5$ and $2\leq k=a<5$. The weight enumerator when
$a=2$ is
\[
W_G(x)=x^5+4x^2+2x+2,
\]
and there are $4$ codewords of minimum weight $3$
in the corresponding elliptic (Goppa) code.
This is a $1$-error correcting code of length $5$
(over $\fff_{3}$).

\subsection{The generator matrix (apr\'es des Goppa)}

\begin{example}
Consider the hyperelliptic curve\footnote{When $p=3$ it is
a model of a modular curve of level $32$ 
(see Table 1). When $p=7$ this 
example arises in the reduction of $X(7)$ 
in characteristic $7$ \cite{E2}.} 
$X$ defined by $y^2=x^p-x$ over the field $\fff_p$ with $p$
elements. It is easy to see that 

\[
C(\fff_p)=\{P_\infty,(0,0),(1,0),...,(p-1,0)\}
\]
has exactly $p+1$ points, including the point at infinity,
$P_\infty$. The automorphism group of this
curve is a two-fold cover of $PSL(2,p)$
(see G\"ob \cite{Go} for the algebraically closed case).

Consider for example the case of $p=7$.
Let $A=mP_\infty$ and $D=P_1+...+P_7$
and let $C$ denote the one-point Goppa code
associated to $X/\fff_7$ and these divisors
$A$, $D$. These codes give rise to MDS codes
in many cases.

When $m=2$, we obtain a 
$[7,2,6]$ code with weight enumerator
$1 + 42x^6 + 6x^7$. 
This code has automorphism group
of order $252$ and permutation group of order $42$.
When $m=4$, we obtain a 
$[7,3,5]$ code with weight enumerator
$1 + 126x^5+ 84x^6 + 132x^7$. 
This code has the same automorphism group
and permutation group.
It has generator matrix in standard form

\[
G=
\left(
\begin{array}{ccccccc}
1 & 0 & 0 & 2 & 5 & 1 & 5\\
0 & 1 & 0 & 1 & 5 & 5 & 2\\
0 & 0 & 1 & 5 & 5 & 2 & 1\\
\end{array}
\right)
\]
and check matrix

\[
H=
\left(
\begin{array}{ccccccc}
5 & 6 & 2 & 1 &0 &0 &0 \\
2 & 2 & 2 & 0 &1 &0 &0 \\
6 & 2 & 5 & 0 &0 &1 &0 \\
2 & 5 & 6 & 0 &0 &0 &1 
\end{array}
\right)
\]
\end{example}

The method used in Goppa's Fermat cubic code
example of \cite{G}, pages 108-109, can be easily
modified to yield analogous
quantities for certain elliptic Goppa codes.

\begin{example}
Let $E$ denote the elliptic curve (of conductor
$N=19$) which we write in
homogeneous coordinates as
\[
y^2z+yz^2=x^3+x^2z+xz^2
\]
Let $\phi(x,y,z)=x^2+y^2+z^2$, let $F$ denote the
projective curve
defined by $\phi(x,y,z)=0$, and let $D$ denote
the divisor obtained by intersecting $E$ and $F$.
By Bezout's theorem, $D$ is
of degree $6$. A basis for ${\cal L}(D)$
is provided by the functions in the set
\[
{\cal B}_D=\{ 1,x^2/\phi(x,y,z),y^2/\phi(x,y,z),
z^2/\phi(x,y,z),xy/\phi(x,y,z),
yz/\phi(x,y,z)\}.
\]
(This is due to the fact that dim ${\cal L}(D)=
{\rm deg}(D)=6$ and the functions $f\in {\cal B}_D$
``obviously'' satisfy $(f)\geq -D$.)
We have
\[
E(\fff_{7})=
\begin{array}{cc}
&\{\ [0, 0, 1], [0, 1, 0], [0, 1, 6], [1, 0, 2],
[1, 0, 4], \\
& [1, 3, 4],
[1, 3, 6], [1, 5, 2], [1, 5, 6]
\ \}
\end{array} ,
\]
which we write as $P_1$, $P_2$, ..., $P_{9}$.
Consider the matrix
\[
G=\left [
\begin {array}{ccccccccc}
0&0&0&1&1&1&1&1&1\\
\noalign{\medskip}0&1&1&0&0&2&2&4&4\\
\noalign{\medskip}1&0&1&4&2&2&1&4&1\\
\noalign{\medskip}0&0&0&0&0&3&3&5&5\\
\noalign{\medskip}0&0&6&0&0&5&4&3&2\\
\noalign{\medskip}0&0&0&2&4&4&6&2&6
\end {array}
\right  ].
\]
The first row of $G$ gives the values of
$x^2/\phi(x,y,z)$ at $\{P_i\ |\ 1\leq i\leq 9\}$.
The other rows
are obtained similarly from the other functions
corresponding to the basis elements of ${\cal L}(D)$:
$y^2/\phi(x,y,z)$, $z^2/\phi(x,y,z)$, $xy/\phi(x,y,z)$,
$yz/\phi(x,y,z)$.
Performing Gauss reduction mod 7 puts this in
canonical form:
\[
G'=
\left [
\begin {array}{ccccccccc}
1&0&0&0&0&0&0&4&4\\
\noalign{\medskip}0&1&0&0&0&0&6&0&6\\
\noalign{\medskip}0&0&1&0&0&0&1&3&4\\
\noalign{\medskip}0&0&0&1&0&0&6&1&6\\
\noalign{\medskip}0&0&0&0&1&0&1&3&5\\
\noalign{\medskip}0&0&0&0&0&1&1&4&4
\end {array}
\right ],
\]
so this code also has minimum distance $3$, hence is
only $1$-error correcting.
The corresponding check matrix is
\[
H=\left [
\begin {array}{ccccccccc}
0 &1 & 2& 1& 2& 2& 1& 0& 0\\
3 &0 & 4& 6& 4& 3& 0& 1& 0\\
3 &1 & 3& 1& 2& 3& 0& 0& 1
\end {array}\right ].
\]

\end{example}

For an example of the generating matrix of a one-point
elliptic code associated to $x^3+y^3=1$ over $\fff_4$ has been
worked out in several places (for example, see Goppa's
book mentioned above, or the books \cite{SS}, \S 3.3,
\cite{P}, \S\S 5.3, 5.4, 5.7, or \cite{M}, \S 5.7.3).

\section{Concluding comments}

        We end this note by making some comments:

        (1) The algebraic geometric relation between the number of points
over a finite field for a variety is related to the Betti numbers. However,
an equivalent notion of genus for higher dimensional varieties is the
``arithmetic genus''. Can one develop a relation between the number of
points over finite fields of a variety and its arithmetic
genus, useful in coding theory?


(2) Can one construct ``good'' codes associated to the higher
dimensional Shimura varieties?

{\bf Acknowledgements}
The first-named author is grateful to 
Pablo Legarraga and Will Traves for useful and
informative discussions. 
The second author wishes to
thank the hospitalities and the financial supports of the Max-Planck
Institut f\"ur Mathematik (Bonn) while he was involved with this
paper. We also thank Amin Shokrollahi for helpful
suggestions.

\end{document}